\documentclass{amsart}
\usepackage[T1]{fontenc}
\usepackage{algorithm,algpseudocode,algorithmicx}
\usepackage[a4paper, left=2.5cm, right=2.5cm, top=3cm, bottom=3cm]{geometry}

\newcommand{\ex}[2]{\ensuremath{\mathrm{ex}(#1, #2)}}
\newcommand{\compoverset}[2]{\ensuremath{K(#2, \overset{#1}{\dots}, #2)}}
\newcommand{\compdots}[2]{\ensuremath{K(#1, \dots, #2)}}
\newcommand{\bigO}[1]{\ensuremath{\mathcal{O}\left(#1\right)}}
\newcommand{\bigOmega}[1]{\ensuremath{\Omega\left(#1\right)}}
\usepackage{parskip}
\usepackage{hyperref}
\usepackage[foot]{amsaddr}

\newtheorem{theorem}[subsection]{Theorem}
\newtheorem{lemma}[subsection]{Lemma}
\newtheorem{remark}[subsection]{Remark}

\title{Finding Partite Hypergraphs Efficiently}

\author[Ferran Espu\~{n}a]{Ferran Espu\~{n}a}
\email{ferran.espuna@gmail.com}

\begin{document}
    \begin{abstract}
    We provide a deterministic polynomial-time algorithm (for fixed $k$) that, for a given $k$-uniform hypergraph $H$ with $n$ vertices and edge density $d$,
    finds a complete $k$-partite subgraph of $H$ with parts of size at least ${c(d, k)(\log n)^{1/(k-1)}}$.
    This generalizes work by Mubayi and Tur\'{a}n on bipartite graphs.
    The value we obtain for the part size matches the order of magnitude guaranteed by the non-constructive proof due to
    Erd\H{o}s and is tight up to a constant factor.
    \end{abstract}

    \maketitle

    \section{Introduction}\label{sec:introduction}

    Hypergraph Tur\'{a}n problems concern how many edges a $k$-uniform hypergraph $H = (V, E)$ with $n$ vertices can have without containing a specific subgraph $G$.
    The maximal such number is known as the \emph{Tur\'{a}n number} $\ex{n}{G}$.
    It is known~\cite{keevash2011hypergraph}
    that $\ex{n}{G}$ is sublinear in $\binom{n}{k}$ if and only if $G$ is $k$-partite, i.e.,
    if its vertex set can be partitioned into $k$ disjoint sets such that each edge contains exactly one vertex from each part.
    K\H{o}v\'{a}ri, S\'{o}s, and Tur\'{a}n~\cite{Kovari1954} (for $k=2$) and
    Erd\H{o}s~\cite{Erdos1964} (for any $k \geq 2$) established that
    \begin{equation} \label{eq:erdos_extremal}
        \ex{n}{\compoverset{k}{t}} \leq \binom{n}{k} \cdot n^{-\frac{1}{t^{k-1}}},
    \end{equation}
    where $\compoverset{k}{t}$ is the hypergraph consisting of $k$ disjoint sets $V_1, \dots, V_k$ of size $t$
    and all hyperedges of the form $\{x_1, \dots, x_k\}$ with $x_i \in V_i$ for $i = 1, \dots, k$.
    This result implies the following.

    \begin{remark} \label{rk:order}
        For all $k \geq 1$ and $0 < d < 1$, there is a constant $c(d, k)$ such that if
        $H$ is a $k$-uniform hypergraph on $n$ vertices with at least $d \binom{n}{k}$ edges,
        and if $t \leq c(d, k) (\log n)^{1/(k-1)}$,
        then $H$ contains $\compoverset{k}{t}$ as a subgraph.
    \end{remark}

    Furthermore, the order of magnitude of $t$ is tight up to a constant factor:
    For some constant $\hat{c}(d, k) > 0$,
    there are $k$-uniform hypergraphs with $n$ vertices and at least $d \binom{n}{k}$ edges
    that do not contain $\compoverset{k}{\hat{t}}$ as a subgraph, as long as
    $\hat{t} \geq \hat{c}(d, k)(\log n)^{1/(k-1)}$.
    This was already noted by Erd\H{o}s~\cite{Erdos1964} and can be proved via the random alteration method~\cite{alon2016probabilistic}.

    Due to the fundamental role of Erd\H{o}s' result,
    it is natural to ask whether a fast search algorithm for a complete $k$-partite subgraph of the size stated in Remark~\ref{rk:order} exists.
    A brute-force search for a $\compoverset{k}{t}$ would involve checking all $\binom{n}{kt}$ vertex subsets,
    which is super-polynomial in the number of vertices $n$ for ${t = \bigOmega{(\log n)^{1/(k-1)}}}$.
    For $k=2$, Mubayi and Tur\'{a}n~\cite{MUBAYI2010174} developed a deterministic polynomial-time algorithm which reaches the stated order of magnitude for the subgraph part size.
    This work extends their approach to the general case of $k$-uniform hypergraphs, reaching analogous results for $k \ge 3$.
    More concretely, we prove the following.

    \begin{theorem} \label{thm:main_theorem}
    There is a deterministic algorithm that, given a $k$-uniform hypergraph $H$ with $n$ vertices and $m=d \binom{n}{k}$ edges, finds a complete balanced $k$-partite subgraph $\compoverset{k}{t}$ in polynomial time, where
    \[
        t = t(n, d, k) = \left\lfloor \left( \frac{\log n}{\log (16/d)}\right)^{\frac{1}{k-1}}\right\rfloor.
    \]
    \end{theorem}
    This value of $t$ matches the asymptotic order of magnitude from Remark~\ref{rk:order}.
    It is worth noting, however, that while our result is asymptotically tight regarding the dependence on $n$,
    the constant factor ${\tilde{c}(d, k) = \log (16/d)^{-\frac{1}{k-1}}}$ obtained by our algorithm is much smaller than that guaranteed by existence proofs.
    For example, $\tilde{c}(d, k)$ remains bounded even as $d \to 1$, whereas the same is not true for the family of constants $c(d, k) = \log(1/d)^{-\frac{1}{k-1}}$
    obtained directly from $\eqref{eq:erdos_extremal}$.
    We focus here on establishing polynomial-time constructibility rather than optimizing this constant to match the best known existential bounds.

    \section{The algorithm}\label{sec:algorithm}

    We present a recursive algorithm, \texttt{FindPartite}, that finds a $\compoverset{k}{t}$ in a given $k$-uniform hypergraph $H$.
    The core idea is to reduce the uniformity of the host hypergraph $H$ from $k$ to $k-1$ in each recursive step.
    The algorithm takes a $k$-uniform hypergraph $H$ with $n$ vertices and $m$ edges as input.
    It first defines the target part size $t$, a small set size $w$,
    and a threshold edge count $s$ for the recursive call, based on the input hypergraph's parameters:
    \begin{align*}
        t &= t(n, d, k) = \left\lfloor \left( \frac{\log n}{\log (16/d)}\right)^{\frac{1}{k-1}}\right\rfloor, \\
        w &= w(n, d, k) = \left\lceil \frac{4 t}{d}\right\rceil \text{, and } \\
        s &= s(n, d, k) = \left\lceil \left( \frac{d}{4}\right)^t \binom{n}{k-1}\right\rceil,
    \end{align*}
    where $d = \frac{m}{\binom{n}{k}}$ is the edge density of $H$.
    The main steps are:
    \begin{enumerate}
        \item \textbf{Base Case ($k=1$):} The edge set of a 1-uniform hypergraph is just a collection of singleton sets of vertices.
        Return the set of all vertices that are ``edges''.

        \item \textbf{Select High-Degree Vertices:} Choose a set $W \subset V$ of $w$ vertices with the highest degrees in~$H$. \label{W}

        \item \textbf{Find a Dense Link Graph:} Iterate through all subsets $T \subset W$ of size $t$.
        For each $T$, consider the set $S$ of all subsets $y \subset V$ of size $k-1$ that form a hyperedge with \emph{every} vertex in $T$. \label{link}

        \item \textbf{Recurse:} As we prove further along using the K\H{o}v\'{a}ri–S\'{o}s–Tur\'{a}n theorem, for at least one choice of $T$,
        the resulting set $S$ is large ($|S| \ge s$). We form a new $(k-1)$-uniform hypergraph $H'=(V, S)$ and make a recursive call: \texttt{FindPartite($H'$, $k-1$)}.

        \item \textbf{Construct Solution:} The recursive call returns $k-1$ parts $V_1, \dots, V_{k-1}$ of size at least $t$.
        Every choice of one vertex from each of these parts forms an edge in $H'$.
        By construction, they form an edge of $H$ with every vertex in $T$.
        Thus, $V_1, \dots, V_{k-1}, T$ span the desired $\compoverset{k}{t}$ subgraph in the original $k$-uniform hypergraph.

    \end{enumerate}

    The pseudocode is given in Algorithm~\ref{alg:kpartite}.
    In section~\ref{sec:correctness}, we formally prove that this algorithm is correct,
    in the sense that it returns a tuple $(V_1, \dots, V_k)$ of disjoint sets $V_i \subset V(H)$
    of size at least $t$ spanning a complete $k$-partite subgraph in $H$.
    In Section~\ref{sec:complexity}, we
    show that it runs in polynomial time in the number of vertices $n$ of the input hypergraph $H$.
    This will conclude the proof of Theorem~\ref{thm:main_theorem}.

    \begin{algorithm}
        \caption{Finding a balanced $k$-partite $k$-uniform subgraph}
        \label{alg:kpartite}
        \begin{algorithmic}[1]
            \Function{FindPartite}{$H, k$}
                \If {$k = 1$}
                    \State \Return $(\{x \colon \{x\} \in E(H)\})$
                \EndIf

                \State $n \gets |V(H)|$, $m \gets |E(H)|$, $d \gets \frac{m}{\binom{n}{k}}$
                \State $t \gets t(n, d, k)$, $w \gets w(n, d, k)$, $s \gets s(n, d, k)$
                \State \textbf{assert} $(t \ge 2)$

                \State $W \gets$ a set of $w$ vertices with highest degree in $H$
                \ForAll{$T \in \binom{W}{t}$}
                    \State $S \gets \{\,y \in \binom{V}{k-1} \colon \forall x \in T, \{x\} \cup y \in E(H)\,\}$
                    \If{$|S| \ge s$}
                        \State $H' \gets (V, S)$  \Comment{$H'$ is a $(k-1)$-uniform hypergraph}
                        \State $(V_1, \dots, V_{k-1}) \gets$ \Call{FindPartite}{$H', k-1$} \label{recurse}
                        \State \Return $(V_1, \dots, V_{k-1}, T)$
                    \EndIf
                \EndFor
            \EndFunction
        \end{algorithmic}
    \end{algorithm}

    \section{Proof of correctness}\label{sec:correctness}

    We now present the proof that
    all the steps in Algorithm~\ref{alg:kpartite} are well-defined and that
    it returns a tuple $(V_1, \dots, V_k)$ of disjoint sets $V_i \subset V(H)$
    of size at least $t$ spanning a complete $k$-partite subgraph in $H$.
    We assume $t \ge 2$ for our estimates to be easier.
    If $t < 2$, we may just return the vertices of any single edge in $H$.

    It is not immediately clear that the set $W$ defined in step~\ref{W} of the algorithm
    is well-defined, as for this it is necessary that $w \leq n$.
    To show this, we first observe that our assumption $t \geq 2$
    implies that $ 1 \geq d \geq \frac{16}{\sqrt{n}}$.
    Suppose, by way of contradiction, that $w > n$.
    Then, we have
    \[
        n \leq w - 1 \leq \frac{4t}{d} \leq \frac{4\log n}{d\log(16/d)} \leq \frac{4 \sqrt{n} \log n }{16 \log(16/d)}.
    \]
    Taking the logarithms to be in base $e$, we note that $\log x \leq \sqrt{x}$ for all positive $x$,
    and that $\log(16/d) > 1$.
    Therefore, we get $n < \frac{n}{4}$, which is a contradiction.

    Next, we prove that in step~\ref{link} of the algorithm
    we indeed find a set $T \in \binom{W}{t}$ such that the associated set
    $S \subset \binom{V}{k-1}$ has size at least $s$.
    That is, Algorithm~\ref{alg:kpartite} reaches line~\ref{recurse} at some point in the for loop.
    For this, consider the bipartite graph $B$ with parts $\binom{V}{k-1}$ and $W$ with edge set
    \[
      \left\{(x, y) \in \binom{V}{k-1} \times W \middle| \, x \cup \{y\} \in E\right\}.
    \]
    The edges of $B$ correspond to the edges containing each vertex in $W$, so there are
    \[
        z = \sum_{y \in W} d_H(y) \geq k \cdot m \cdot \frac{w}{n} = \frac{k \cdot w \cdot d \cdot \binom{n}{k}}{n} = w \cdot d \cdot \binom{n - 1}{k-1}
    \]
    of them, where the inequality follows from the fact that we have picked a set of $w$ vertices with highest degree in $H$.
    The existence of a set $T \subset W$ as desired is equivalent to there being $T \subset W$ of size $t$ and a set $S \subset \binom{V}{k-1}$ of size $s$
    such that the induced bipartite subgraph $B[S, T]$ is complete.
    To prove that this is the case,
    we use the following version~\cite{Hylten1958} of the K\H{o}v\'{a}ri–S\'{o}s–Tur\'{a}n theorem~\cite{Kovari1954}.

    \begin{lemma}\label{thm:kst}
        Let $u, w, s, t$ be positive integers with $u \geq s$, $w \geq t$, and let $B$ be a bipartite graph with parts $W$ and $U$ such that
        $|U| = u, |W| = w$.
        If $B$ has more than \[(s - 1)^{1 / t}(w - t + 1)u^{1 - 1 / t} + (t - 1)u\] edges, then there are
        $T \subset W$ of size $t$ and $S \subset U$ of size $s$ such that the induced bipartite subgraph $B[S, T]$ is complete.
    \end{lemma}

    We apply the lemma with $u = \binom{n}{k-1}$.
    It is clear from the definitions that our parameters satisfy the requirements $u \geq s$ and $w \geq t$.
    Suppose, by way of contradiction, that
    \[
        w \cdot d \cdot \binom{n - 1}{k-1} \leq z \leq (s - 1)^{1 / t}(w - t + 1)\binom{n}{k-1}^{1 - 1 / t} + (t - 1)\binom{n}{k-1}.
    \]
    Dividing by $\binom{n}{k-1}$ then shows that
    \[
        \frac{1}{2} \cdot w \cdot d
        \leq \left( 1 - \frac{k-1}{n}\right) \cdot w \cdot d
        = w \cdot d \cdot \frac{\binom{n-1}{k-1}}{\binom{n}{k-1}}
        < w \left( \frac{s-1}{\binom{n}{k-1} }\right)^{1 / t} + (t - 1),
    \]
    where the first inequality follows from $n \geq 16^{2^{k-1}} > 2(k-1)$, which follows from $t \geq 2$ and $d \leq 1$.
    Finally, since $t \leq \frac{w \cdot d}{4}$ by the definition of $w$, we obtain
    \[
        \left( \frac{d}{4}\right)^t \binom{n}{k-1} < s-1,
    \]
    in contradiction to the definition of $s$.
    So far, we have shown that the sets defined in Algorithm~\ref{alg:kpartite}
    are well-defined and that the recursive call in line~\ref{recurse} is reached.
    We are now ready to prove that the algorithm returns a $\compoverset{k}{t}$
    by examining what happens in the recursive call.
    More precisely, we show the following.

    \begin{theorem}
        For $k \geq 2$, if $t \geq 2$, Algorithm~\ref{alg:kpartite} returns a tuple $(V_1, \dots, V_k)$ of disjoint sets $V_i \subset V(H)$ such that
        $|V_i| \geq t$ and $H[V_1, \dots, V_k]$ is complete.

        \begin{proof}
            We proceed by induction on $k$.
            For $k=2$, the recursive call returns the common neighborhood $V_1$ of the vertices in $T$,
            which is obviously disjoint from $T$, so it only remains to check that $|V_1| \geq t$.
            Now, since by construction $|V_1| = |S| \geq s$, it is enough that
            \[
                s
                = \left\lceil \left( \frac{d}{4}\right)^t \cdot n\right\rceil
                \geq \left( \frac{d}{4}\right)^{\frac{\log n}{\log(16/d)}}  \cdot n
                = \left( 4 \cdot \frac{d}{16}\right)^{\frac{\log n}{\log(16/d)}} \cdot n
                = 4^{\frac{\log n}{\log(16/d)}} \cdot \frac{1}{n} \cdot n
                \geq 4^t
                > t.
            \]

            For $k \geq 3$, we assume the inductive hypothesis holds for $k-1$.
            Let $d'$ be the edge density of $H'$ (which is the obtained $(k-1)$-uniform hypergraph)
            and $t'$ be $t(n, d', k-1)$.
            As long as $t' \geq 2$, the recursive call returns
            a tuple $(V_1, \dots, V_{k-1})$ of disjoint sets $V_i \subset V(H)$ such that
            $|V_i| \geq t'$ and $H'[V_1, \dots, V_{k-1}]$ is complete.

            We claim that $t' \geq t$.
            This implies that $t' \geq 2$ so we get to apply the inductive hypothesis to $H'$.
            Furthermore, the sets $V_i$ that we obtain when applying the algorithm to $H'$
            are of size at least $t' \geq t$.
            By construction of $S$, all the edges in $H'$ are disjoint from $T$,
            hence the sets $V_i$ are also disjoint from $T$.
            This means that the sets $V_1, \dots, V_{k-1}, T$ are disjoint.
            In addition, for all $(x_1, \dots x_{k-1}, y) \in V_1 \times \dots \times V_{k-1} \times T$,
            we have that $\{x_1, \dots, x_{k-1}\} \in S$ so $\{x_1, \dots, x_{k-1}, y\} \in E(H)$.
            Equivalently, $H[V_1, \dots, V_{k-1}, T]$ is complete, finishing the proof.

            Let us now prove the claim that $t' \geq t$.
            By the definition of $s$, we have $d' \geq \left( \frac{d}{4}\right)^t$.
            Therefore,
            \[
                t' \geq
                \left\lfloor \left(  \frac{\log n}{\log \left(\frac{16}{(d/4)^t}\right)}\right)^
                {\frac{1}{k-2}}\right\rfloor =
                \left\lfloor \left(  \frac{\log n}{\log 16 + t \log (4/d)}\right)^{\frac{1}{k-2}}\right\rfloor.
            \]
            Then, we substitute the definition of $t$, where removing the floor function
            maintains the inequality because the right side is decreasing in $t$:
            \begin{equation*}
                t' \geq
                \left\lfloor \left(  \frac{\log n}
                {\log 16 + \left(  \frac{\log n}{\log (16/d)}\right)^{\frac{1}{k-1}}  \log (4/d)}\right)^
                {\frac{1}{k-2}}\right\rfloor
                =
                \left\lfloor \left(  \frac{(\log n)^{\left(1-\frac{1}{k-1}\right)}}
                {\frac{\log 16}{(\log n)^{\frac{1}{k-1}}} + \frac{\log (4/d)}{\log (16/d)^{\frac{1}{k-1}}} }
               \right)^{\frac{1}{k-2}}\right\rfloor.
            \end{equation*}
            We claim that we can bound the denominator by showing that
            \begin{equation} \label{eq:denominator_bound}
                \frac{\log 16}{(\log n)^{\frac{1}{k-1}}} + \frac{\log (4/d)}{\log (16/d)^{\frac{1}{k-1}}}
                \leq \left( \log (16/d)\right)^{\left( 1 - \frac{1}{k-1}\right)}.
            \end{equation}
            Then, the expression simplifies to
            \[
                t'
                \geq \left\lfloor \left(  \frac{(\log n)^{\left(1-\frac{1}{k-1}\right)}}
                {\left( \log (16/d)\right)^{\left( 1 - \frac{1}{k-1}\right)}}
                \right)^{\frac{1}{k-2}}\right\rfloor
                = \left\lfloor \left(  \frac{\log n}
                {\log (16/d)}
                \right)^{\frac{1}{k-2}\left( 1 - \frac{1}{k-1}\right)}\right\rfloor
                = \left\lfloor \left(  \frac{\log n}{\log (16/d)}\right)^{\frac{1}{k-1}}\right\rfloor
                = t,
            \]
            as desired.
            Let us prove Inequality~\eqref{eq:denominator_bound}.
            Suppose, by way of contradiction, that it does not hold.
            We can rewrite
            \[
                (\log (16/d))^{\left( 1 - \frac{1}{k-1}\right)}
                = \frac{\log (16/d)}{\log (16/d)^{\frac{1}{k-1}}}
            \]
            and rearrange the inequality to obtain
            \[
                \frac{\log 16}{(\log n)^{\frac{1}{k-1}}}
                > \frac{\log (16/d) - \log (4/d)}{\log (16/d)^{\frac{1}{k-1}}}
                = \frac{\log 4}{(\log (16/d))^{\frac{1}{k-1}}}.
            \]
            This implies that
            \[
                t
                \leq \left( \frac{\log n}{\log (16/d)}\right)^{\frac{1}{k-1}}
                < \frac{\log 16}{\log 4} = 2,
            \]
            which contradicts the assumption that $t \geq 2$.
        \end{proof}
    \end{theorem}

    \section{Proof of Polynomial Complexity}\label{sec:complexity}

    We now analyze the computational complexity of Algorithm~\ref{alg:kpartite} to show that it runs in time polynomial in $n$ for any fixed uniformity $k$.
    We consider the input hypergraph $H$ to be given as an adjacency $k$-dimensional tensor.
    If the hypergraph is provided sparsely (e.g., as a list of $m$ edges), it can be preprocessed in polynomial time by reindexing the vertices and removing isolated ones.
    This reduces the effective number of vertices to $n'$ such that $n' \le k \cdot m$.
    Consequently, the tensor initialization time and space become $\bigO{(km)^k}$, which remains polynomial in the sparse input size for fixed $k$.
    Since we will prove the algorithm is polynomial-time in this new size, the resulting algorithm is also polynomial in the sparse size.

    We perform the analysis under the standard Word RAM model with word size $\Omega(\log n)$.
    This allows for basic arithmetic operations and array indexing on values up to $n^k$ to be performed in $\bigO{1}$ time.

    Let $T_k(n)$ denote the worst-case running time of the function \texttt{FindPartite} when called on a $k$-uniform hypergraph with $n$ vertices.
    The algorithm's structure gives a recurrence relation for $T_k(n)$.
    We first analyze the cost of the operations within a single call for a fixed $k$, excluding the recursive step.

    Querying whether a set of $k$ vertices forms an edge in $H$ can be done in constant time.
    Therefore, calculating the number of edges $m$ can be done in $\bigO{n^k}$ time by iterating through all $\binom{n}{k}$ possible edges.
    The parameters $t, w$, and $s$ can then be calculated in polylogarithmic time: If one wants to avoid the arithmetic on real
    numbers, the values can be obtained via binary search on equations obtained by rearranging the definitions of $t, w$, and $s$.

    The set $W$ of $w$ vertices with highest degrees can be constructed in time $\bigO{n^k}$, for instance, by
    creating an array with the degree of each vertex (in time  $\bigO{n \cdot n^{k-1}} = \bigO{n^k}$),
    sorting it in descending order (in time $\bigO{n \log n}$),
    and taking the first $w$ elements.

    The subsets of $t$ elements of $W$ can be iterated through in time $\bigO{\binom{w}{t}}$~\cite{reingold1977combinatorial}.
    Similarly to the argument by Mubayi and Tur\'{a}n~\cite{MUBAYI2010174}, we can bound this quantity by a polynomial in $n$.
    Indeed,
    \[
        \binom{w}{t} \leq \left(\frac{ew}{t}\right)^t
        \leq \left(\frac{e(4t/d + 1)}{t}\right)^t
        \leq \left(\frac{4e}{d} + \frac{e}{2}\right)^t
        \leq \left(\frac{4.5 \cdot e}{d}\right)^t < e^{t(2.6 + \log (1/d))}
        = e^{2.6 \cdot t} e^{\log (1/d) \cdot t}.
    \]
    For the first term, we use the fact that $t < (\log n)^{1/(k-1)} \leq \log n$
    and so ${e^{2.6 \cdot t} < e^{2.6 \cdot \log n} = n^{2.6}}$.
    For the second term, we use that $t < ((\log n)/\log(1/d))^{1/(k-1)} \leq (\log n)/\log(1/d)$,
    so $e^{\log (1/d) \cdot t} < n$.
    All in all, we have
    \[
        \binom{w}{t} < n^{2.6} n = n^{3.6}.
    \]
    In each iteration of the loop, the set $S$ can be constructed in the following way.
    We initialize a boolean adjacency tensor $A$ of size $n^{k-1}$,
    with all entries set to \texttt{true}.
    We iterate through all $x \in \binom{[n]}{k-1}$ and $y \in T$
    and set the entry corresponding to $x$ to \texttt{false} if $x \cup \{y\}$ is not an edge in $H$.
    All in all, this takes $\bigO{tn^{k-1}} = \bigO{n^k}$ steps, and then counting the number of \texttt{true} entries in the array takes $\bigO{n^{k-1}}$ time.
    Therefore, the for loop (without the recursive call) takes $\bigO{\binom{w}{t} n^{k}} = \bigO{n^{k+3.6}}$ time.

    Finally, when the condition $|S| \geq s$ is satisfied, the recursive call to \texttt{FindPartite} is made.
    We can pass the array $A$ to the recursive call directly, and the recursive call takes time $T_{k-1}(n)$.
    Putting everything together, we have the recurrence relation $T_k(n) = T_{k-1}(n) + \bigO{n^{k+3.6}}$.
    This, together with the base case $T_1(n) = \bigO{n}$, gives us $T_k(n) = \bigO{n^{k+3.6}}$.
    Similar analysis shows that if the more restrictive Turing Machine computational model is used,
    the algorithm runs in time $\bigO{n^{2k + 3.6} p_k(\log(n))}$ where $p_k$ is some polynomial depending on $k$.

    \section{Concluding Remarks and Future Work}\label{sec:conclusion-and-future-work}

    We have presented a deterministic,
    polynomial-time algorithm to find a large complete balanced $k$-partite subgraph in any sufficiently dense $k$-uniform hypergraph.
    This provides a constructive counterpart to a classical existence result by Erd\H{o}s in extremal hypergraph theory.

    Several avenues for future research remain open.
    \begin{itemize}
        \item \textbf{General Blow-ups:} Our algorithm finds a blow-up of a single edge, $\compoverset{k}{t}$.
        Can this framework be adapted to find a $t$-blowup of an arbitrary fixed $k$-uniform hypergraph $G$ when its Tur\'{a}n density is exceeded by a positive constant?
        For example, for $k=2$, there are existence results for $t = \bigOmega{\log n}$~\cite{bollobas1973structure}, but directly adapting Algorithm~\ref{alg:kpartite} would only
        yield $t = \bigO{(\log n)^{1/(|V(G)| - 1)}}$.
        \item \textbf{Unbalanced k-Partite Graphs:} The algorithm could be modified to search for unbalanced complete partite graphs $\compdots{t_1}{t_k}$, where the part sizes may grow at different rates.
        \item \textbf{Optimality:} The bounds on $t$ are asymptotically tight, but the constants can be improved significantly with a more refined analysis.
    \end{itemize}

    \subsection*{Acknowledgements}

    The ideas in this work stem from the author's master's thesis of the same name at
    Universitat Polit\`{e}cnica de Catalunya,
    under the supervision of Richard Lang.
    The author thanks him for suggesting the problem and for his guidance and support throughout the project.

    \bibliography{journal}
    \bibliographystyle{plain}

\end{document}